\documentclass[12pt]{article}

\usepackage{amsmath, epsfig, cite}
\usepackage{amssymb}
\usepackage{amsfonts}
\usepackage{latexsym}
\usepackage{float}
\usepackage{color}
\usepackage{mathrsfs}

\newtheorem{thm}{Theorem}[section]

\newtheorem{cor}[thm]{Corollary}

\newtheorem{lem}[thm]{Lemma}

\newcommand{\pf}{\noindent{\it Proof} }

\numberwithin{equation}{section}

\newcommand{\qed}{{\hfill$\square$}\medskip}

\begin{document}
	%\linenumbers
	\begin{center}
		{\large\bf  A Generalized Supercongruence  of  Z.-W. Sun   }
	\end{center}
	\vskip 2mm \centerline{Wei-Wei Qi}
	
	\begin{center}
		{\footnotesize MOE-LCSM, School of Mathematics and Statistics, Hunan Normal University, Hunan 410081, P.R. China\\[5pt]
			{\tt wwqi2022@foxmail.com} \\[10pt]
		}
	\end{center}
	
	\vskip 0.7cm \noindent{\bf Abstract.} In this paper, we  employ the Wilf-Zeilberger (WZ) method  to prove a supercongruence conjecture  posed by Z.-W. Sun: for any prime $p$,	
	\begin{align*}
		\sum_{k=0}^{\frac{p-3}{2}}\frac{92k^2+61k+9}{(2k+1)64^k}{2k \choose k}{3k \choose k}{4k \choose 2k}\equiv 6p+16p^2\left(\frac{-1}{p}\right) \pmod{p^3},
	\end{align*}	 
where $\left(\frac{\cdot}{p}\right)$ denotes the Legendre symbol. Our proof relies on  combinatorial identities and symbolic summation techniques.

	\vskip 2mm
	\noindent{\it MR Subject Classifications}: 11A07, 05A10, 05A19.	
	
\section{Introduction} 
	
In $1914$, Ramanujan \cite{www-1-5} presented a number of infinite series for $1/\pi$  in the form
\begin{align*}
 \sum_{k=0}^{\infty}(ak+b)\frac{t_k}{m^k}=\frac{c}{\pi},
\end{align*}		
which $a,b,c,m$ are constants and $t_k$ denotes a product of the binomial coefficients ${2k\choose k}$, ${3k\choose k}$, ${4k \choose 2k}$. In $1997$, Van Hamme \cite{www-1-6} investigated the partial sums of Ramanujan-type series for $1/\pi$, and put forward thirteen conjectured supercongruences labeled as $(A.2)-(M.2)$. Up to now, all of these thirteen supercongruences has been independently confirmed by  various distinct methods. In $1990$, H. S. Wilf and D. Zeilberger \cite{www-1-7} introduced the WZ method, which can  discover and prove hypergeometric identities. In $1995$, Gessel developed a systematic approach to construct WZ pairs. Recently, Feng and Hou (\cite{www-2}, \cite{www-2-0}) established several congruences related to truncated Ramanujan-type series by applying the WZ method. For more applications of WZ method, readers may refer to \cite{www-1-9}, \cite{www-1-10}, \cite{www-1-11}  and other relevant literature.

 In $2010$, Sun and Tauraso \cite{www-1-12} showed that for prime $p$, $d\in \{0,\cdots, p^a\}$ with $a\in \mathbb{Z}^+$. Let $m\in \mathbb{Z}$ with $p\nmid m$, then
\begin{align*}
	\sum_{k=0}^{p^a-1}\frac{{2k\choose k+d}}{m^k} \equiv u_{p^a-d}(m-2) \pmod{p},
\end{align*}	
where $u_{-1}(x)=xu_0(x)-u_1(x)=-1$. For a fixed integer $x$, the sequence $\{u_n(x)\}_{n\in \mathbb{N}}$ is linear recurrence of integer. Namely, for any $n\in\mathbb{N}$, the equality $u_n(-x)=(-1)^{n-1}u_n(x)$ holds. In $2011$, Sun \cite{www-1-13} established that
\begin{align*}
	\sum_{k=0}^{\frac{p-1}{2}}\frac{{2k\choose k}}{k} \equiv -\left(\frac{-1}{p}\right)\frac{8}{3}E_{p-3} \pmod{p^2},
\end{align*}
where $\left(\frac{\cdot}{p}\right)$ denotes the Legendre symbol  and  $E_{p-3}$ stands for the
  ($p-3$)-Euler numbers. And the  Euler numbers is given by
\begin{align*}
	\frac{2}{e^x+e^{-x}}=\sum_{n=0}^{\infty}E_n\frac{x^n}{n!}.
\end{align*}
Congruences involving central  binomial coefficients are an interesting project, which have been widely studied by many authors. For more details congruences involving central binomial coefficients, one may consult (\cite{www-1-13}, \cite{www-1-14}, \cite{www-1-15}, \cite{www-1-16} and so on). 

 Furthermore, Sun \cite{www-1} summarized numerous conjectures on congruences. For example, let $p$ be an odd prime \cite[Conjecture 4.15]{www-1}, then	
\begin{align}
	\sum_{k=0}^{\frac{p-3}{2}}&\frac{92k^2+61k+9}{(2k+1)64^k}{2k \choose k}{3k \choose k}{4k \choose 2k}
	\equiv 6p+16p^2\left(\frac{-1}{p}\right) \pmod{p^3}. \label{01}
\end{align}		

The first goal of this paper is to extend the congruence conjecture \eqref{01} to a supercongruence modulo $p^4$. The Fermat quotient of an integer $a$ with respect to an odd prime $p$ is gives by
\begin{align*}
	q_p(a)=\frac{a^{p-1}-1}{p}.
\end{align*}

	\begin{thm}
		For any  prime  integer $p>2$. Then
\begin{align}
	\begin{aligned}
	\sum_{k=0}^{\frac{p-3}{2}}&\frac{92k^2+61k+9}{(2k+1)64^k}{2k \choose k}{3k \choose k}{4k \choose 2k}\\
	&\equiv 6p+\left(16p^2+52p^3-48p^3q_p(2)\right)\left(\frac{-1}{p}\right)-12p^3\left(\frac{-1}{p}\right)E_{p-3} \pmod{p^4}. \label{03}
\end{aligned}	
\end{align}	
	\end{thm}
{\it Remark}: We immediately obtain \eqref{01} from \eqref{03}.

Feng and Hou \cite[Theorem 1.4]{www-2} employed the  WZ method to establish the following congruence: for any prime $p$,
	\begin{align*}
	\sum_{k=0}^{p-1}\frac{92k^2+61k+9}{(2k+1)64^k}{2k \choose k}{3k \choose k}{4k \choose 2k}\equiv 0 \pmod{p}. 
\end{align*}
They remark that Sun \cite[Conjecture 4.15]{www-1} conjectured this congruence modulo $p^3$. Here,  our second aim is to investigate this consequence.	
\begin{thm}
Let $p$ be an odd prime. Then
\begin{align}
	\sum_{k=0}^{p-1}\frac{92k^2+61k+9}{(2k+1)64^k}{2k \choose k}{3k \choose k}{4k \choose 2k}\equiv 9p \pmod{p^3}. \label{04}
\end{align}
\end{thm}
Inspired by the proof of Theorem $1.2$, we deduce the follows consequences:
\begin{cor}
	Let $p$ be an odd prime. Then
	\begin{align}
		\sum_{k=1}^{\frac{p-1}{2}}\frac{{2k \choose k}}{2^k}\equiv  -\frac{p}{2}\left(\frac{-1}{p}\right)\sum_{k=1}^{\frac{p-1}{2}}\frac{2^k}{k}+\left(\frac{-1}{p}\right)-1 \pmod{p^2},  \label{06}	
	\end{align}
and
\begin{align}
	\sum_{k=1}^{\frac{p-1}{2}}\frac{{2k \choose k}}{2^k}\equiv  \left(\frac{-1}{p}\right)-1 \pmod{p}.	\label{05}
\end{align}
\end{cor}	
Obviously, congruence \eqref{05} follows from congruence \eqref{06}.	
	
The structure of this paper is organized as follows. In Section $2$, We present some necessary Lemmas. We are going to prove Theorem $1.1$, Theorem $1.2$ and Corollary $1.3$ in Section $3$,  $4$ and  $5$, respectively.	
		
\section{Some Lemmas}
Let $H_n(r)$ denote the $n$-th generalized harmonic number of order $r$, which is defined as
\begin{align*}
	H_n(r)=\sum_{k=1}^{n}\frac{1}{k^r},
\end{align*}
with the convention that $H_n=H_n(1)$. Additionally, we adopt the following notation:	
\begin{align*}
 H_n(-r)=\sum_{k=1}^{n}\frac{(-1)^k}{k^r} \quad and \quad H(1,1;n)=\sum_{1\leq i<j\leq n}\frac{1}{ij}.
\end{align*}

	\begin{lem} For any prime $p>2$, we have	

	\begin{align}
		\sum_{k=1}^{\frac{p-1}{2}}\frac{(-1)^kH_k}{k}&\equiv \frac{1}{2}q_p(2)^2+\left(\frac{-1}{p}\right)E_{p-3} \pmod{p},  \label{0l-8}\\
		\sum_{k \underset{4|k}=1}^{p-1}\frac{H_k}{k}&\equiv \frac{5}{16}q_p(2)^2 \pmod{p},  \label{0l-9}\\
		\sum_{k=1}^{\frac{p-1}{2}}\frac{H_{2k}}{k}&\equiv	q_p(2)^2 \pmod{p}.  \label{0l-10}
	\end{align}	
	
	\end{lem}	
	\pf. The  congruences \eqref{0l-8}, \eqref{0l-9} from \cite[(3.12), (3.16)]{www-1-2} and \eqref{0l-10} from \cite[(2.4)]{www-1-3}.	\qed

	\begin{lem} 
		Let $p>2$ be an odd prime. Then
		\begin{align}
			\begin{aligned}
				\sum_{k=1}^{\frac{p-7}{2}} \frac{(-1)^k}{(2k+1)(2k+3)(2k+5)}
				\equiv  \left(\frac{-1}{p}\right)\left(\frac{q_p(2)}{4}-\frac{p}{8}q_p(2)^2-\frac{11}{128}p-\frac{5}{32}\right)-\frac{2}{5}    \pmod{p^2}, \label{0l-1-0}
			\end{aligned}	
		\end{align}	
	\begin{align}
		\begin{aligned}
			\sum_{k=1}^{\frac{p-7}{2}}& \frac{(-1)^kH_k}{(2k+1)(2k+3)(2k+5)}\\
			&\equiv  \left(\frac{-1}{p}\right)\left(\frac{35}{48}q_p(2)-\frac{1}{2}q_p(2)^2-\frac{77}{72}\right)+\frac{1}{3}q_p(2)+\frac{1}{4}E_{p-3}-\frac{29}{36}    \pmod{p}, \label{0l-1}
		\end{aligned}	
	\end{align}	
	and
	\begin{align}
		\begin{aligned}
			\sum_{k=1}^{\frac{p-7}{2}} \frac{(-1)^kH_{2k}}{(2k+1)(2k+3)(2k+5)}
			\equiv  \left(\frac{-1}{p}\right)\left(\frac{19}{48}q_p(2)-\frac{1}{16}q_p(2)^2-\frac{13}{18}\right)+\frac{1}{6}q_p(2)-\frac{11}{18}    \pmod{p}.  \label{0l-2}
		\end{aligned}	
	\end{align}		
	\end{lem}	
\pf.  Recall the classical Wolstenholme congruence \cite{www-1-0} says that
\begin{align*}
		H_{p-1}\equiv 0   \pmod{p^2}, 
\end{align*}
and a result of Lehmer \cite{www-1-1} asserts that
\begin{align}
	\begin{aligned}
		H_{(p-1)/2}\equiv -2q_p(2)+q_p(2)^2p   \pmod{p^2}. \label{0l-4}
	\end{aligned}	
\end{align}

Since 
\begin{align*}
	\begin{aligned}
		\sum_{k=1}^{p-1}\frac{(-1)^k}{k}=	\sum_{k=1}^{\frac{p-1}{2}}\frac{(-1)^k}{k}+	\sum_{k=1}^{\frac{p-1}{2}}\frac{(-1)^{p-k}}{p-k}\equiv 2	\sum_{k=1}^{\frac{p-1}{2}}\frac{(-1)^k}{k}   \pmod{p^2}, 
	\end{aligned}	
\end{align*}
hence combining  the above congruence and fact $H_{p-1}(-1)=H_{(p-1)/2}-H_{p-1}$, we have
\begin{align}
	\begin{aligned}
		H_{(p-1)/2}(-1)\equiv -q_p(2)   \pmod{p}.  \label{0l-6}
	\end{aligned}	
\end{align}

In view of \cite[Lemma 2.2]{www-1-4}, we can get
\begin{align}
	\begin{aligned}
		\sum_{k=1}^{\frac{p-1}{2}}\frac{1}{k^2}\equiv \frac{1}{2}\sum_{k=1}^{p-1}\frac{1}{k^2}\equiv 0 \pmod{p}.  \label{03-3-1}
	\end{aligned}	
\end{align}
In \cite{www-1-2}, Sun also showed that
\begin{align*}
	\begin{aligned}	
		\sum_{k=1}^{\lfloor \frac{p}{4}\rfloor}\frac{1}{k^2} &\equiv 4\left(\frac{-1}{p}\right)E_{p-3} \pmod{p},\\
		\sum_{k=1}^{\lfloor \frac{p}{4}\rfloor}\frac{1}{k} &\equiv -3q_p(2)+\frac{3}{2}pq_p(2)^2 -\left(\frac{-1}{p}\right)pE_{p-3} \pmod{p^2}. 
	\end{aligned}	
\end{align*}
This, together with   \eqref{0l-4} and  \eqref{03-3-1}  yields
\begin{align}
	\begin{aligned}
		\sum_{k=1}^{\frac{p-3}{2}}\frac{(-1)^k}{1+2k}&=(-1)^{\frac{p-1}{2}}\left(2\sum_{k=1}^{\lfloor \frac{p-1}{4}\rfloor}\frac{1}{p-4k}-\sum_{k=1}^{\frac{p-1}{2}}\frac{1}{p-2k} \right)-1\\
		&
		\equiv (-1)^{\frac{p-1}{2}}\frac{q_p(2)}{2}-(-1)^{\frac{p-1}{2}}\frac{p}{4}q_p(2)^2-1   \pmod{p^2}.   \label{0l-7}
	\end{aligned}	
\end{align}	
By \eqref{0l-7},  we obtain
\begin{align*}
	\begin{aligned}
		\sum_{k=1}^{\frac{p-7}{2}}& \frac{(-1)^k}{(2k+1)(2k+3)(2k+5)}\\
		&=\frac{1}{2}\sum_{k=1}^{\frac{p-3}{2}}\frac{(-1)^k}{1+2k}-\frac{(-1)^{\frac{p-1}{2}}(p-5)}{4(p-2)(p-4)}\\
		&\equiv  \left(\frac{-1}{p}\right)\left(\frac{q_p(2)}{4}-\frac{p}{8}q_p(2)^2-\frac{11}{128}p-\frac{5}{32}\right)-\frac{2}{5}    \pmod{p^2},
	\end{aligned}	
\end{align*}	
as claimed result.

Moreover, with the help of \eqref{0l-8}--\eqref{0l-10}, \eqref{0l-4} and \eqref{0l-6}  gives
\begin{align}
	\begin{aligned}
		\sum_{k=1}^{\frac{p-3}{2}}\frac{(-1)^kH_k}{1+2k}&=-(-1)^{\frac{p-1}{2}}
		\sum_{k=1}^{\frac{p-1}{2}}\frac{(-1)^kH_{\frac{p-1}{2}}-k}{p-2k}\\
		&\equiv (-1)^{\frac{p-1}{2}}\left(\frac{1}{2}\sum_{k=1}^{\frac{p-1}{2}}\frac{(-1)^kH_k}{k}-\sum_{k=1}^{\frac{p-1}{2}}\frac{(-1)^kH_{2k}}{k}-H_{(p-1)/2}(-1)H_{(p-1)/2}\right)\\
		&= (-1)^{\frac{p-1}{2}}\left(\frac{1}{2}\sum_{k=1}^{\frac{p-1}{2}}\frac{(-1)^kH_k}{k}+\sum_{k=1}^{\frac{p-1}{2}}\frac{H_{2k}}{k}-\sum_{k=1}^{\lfloor\frac{p-1}{4}\rfloor}\frac{H_{4k}}{k}
		-H_{(p-1)/2}(-1)H_{(p-1)/2}\right)\\
		&\equiv -(-1)^{\frac{p-1}{2}}q_p(2)^2+\frac{1}{2}E_{p-3}   \pmod{p}, \label{0l-11}
	\end{aligned}	
\end{align}		
where  we used congruence ($H_{(p-1)/2-k}\equiv H{(p-1)}+2H_{2k}-H_k\quad \pmod{p}$) in the second step.	
	
Next, using \eqref{0l-4}, \eqref{0l-6}, \eqref{0l-7} and \eqref{0l-11}, we get	
	\begin{align*}
	\begin{aligned}
		\sum_{k=1}^{\frac{p-7}{2}}& \frac{(-1)^kH_k}{(2k+1)(2k+3)(2k+5)}\\
		&=\sum_{k=1}^{\frac{p-3}{2}}\frac{(-1)^k}{1+2k}\left(\frac{H_k}{2}+\frac{5}{6}\right)-\frac{1}{3}H_{\frac{p-1}{2}}(-1)+\frac{(-1)^{\frac{p-1}{2}}(p-5)}{4(p-2)(p-4)}H_{\frac{p-1}{2}}   \\
	&+(-1)^{\frac{p-1}{2}}\frac{9p^3-88p^2+291p-308}{12(p-1)(p-2)(p-3)(p-4)}+\frac{1}{36}\\
		&\equiv  \left(\frac{-1}{p}\right)\left(\frac{35}{48}q_p(2)-\frac{1}{2}q_p(2)^2-\frac{77}{72}\right)+\frac{1}{3}q_p(2)+\frac{1}{4}E_{p-3}-\frac{29}{36}    \pmod{p}, 
	\end{aligned}	
\end{align*}		
as desired consequence.	

Similarly, utilizing ($H_{p-2k-1}(1)\equiv H_{2k} \quad \pmod{p}$), \eqref{0l-9}, \eqref{0l-10}, \eqref{0l-6} and  \eqref{0l-7}  we obtain
\begin{align}
	\begin{aligned}	
		\sum_{k=1}^{\frac{p-3}{2}}\frac{(-1)^kH_{2k-1}}{1+2k}&=(-1)^{\frac{p-1}{2}}\sum_{k=1}^{\frac{p-1}{2}}\frac{(-1)^kH_{p-2k-1}}{p-2k}-\sum_{k=1}^{\frac{p-1}{2}}\frac{(-1)^k}{2k}+\sum_{k=1}^{\frac{p-3}{2}}\frac{(-1)^k}{1+2k}+\frac{(-1)^{\frac{p-1}{2}}}{p-1}\\
		&\equiv -\frac{(-1)^{\frac{p-1}{2}}}{8}q_p(2)^2+\frac{1+(-1)^{\frac{p-1}{2}}}{2}q_p(2)-(-1)^{\frac{p-1}{2}}-1 \pmod{p}.  \label{0l-12}
	\end{aligned}	
\end{align}
Clearly, we have
	\begin{align*}
	\begin{aligned}
		\sum_{k=1}^{\frac{p-7}{2}}& \frac{(-1)^kH_{2k}}{(2k+1)(2k+3)(2k+5)}\\
		&=\frac{1}{8}\sum_{k=1}^{\frac{p-7}{2}} \frac{(-1)^kH_{2k}}{2k+1}-\frac{1}{4}\sum_{k=1}^{\frac{p-7}{2}} \frac{(-1)^kH_{2k}}{2k+3}+\frac{1}{8}\sum_{k=1}^{\frac{p-7}{2}} \frac{(-1)^kH_{2k}}{2k+5} \\
		&\equiv \frac{1}{2}\sum_{k=1}^{\frac{p-3}{2}} \frac{(-1)^kH_{2k-1}}{2k+1} +\frac{7}{24}\sum_{k=1}^{\frac{p-3}{2}} \frac{(-1)^k}{2k+1}+\frac{H_{\frac{p-1}{2}}(-1)}{12}-\frac{2}{9}(-1)^{\frac{p-1}{2}}+\frac{13}{72}     \pmod{p}.  
	\end{aligned}	
\end{align*}		
Then, substituting \eqref{0l-6}, \eqref{0l-7} and \eqref{0l-12} into above congruence and simplifying, we arrive at \eqref{0l-2}. Hence the proof is complete. \qed
	
	\begin{lem} 
	For any positive integer $n$, we have
	\begin{align}
		\sum_{k=1}^{n}\frac{{n\choose k}}{1+2k}=\frac{4^n}{(1+2k){2n \choose n}}\sum_{k=1}^{n}\frac{{2k\choose k}}{2^k}+\frac{4^n}{(1+2n){2n\choose n}}-1, \label{h-1}
	\end{align}		
and
		\begin{align}
		\sum_{k=1}^{n}\frac{{n\choose k}}{k}=-\sum_{k=1}^{n}\frac{1}{k}+\sum_{k=1}^{n}\frac{2^k}{k}. \label{h-2}
	\end{align}
\end{lem}
\pf. Actually, these identities can be derived by purely combinatorial means. However, using Schneider's \cite{www-3} computer algebra package \textbf{Sigma}, we more quickly find that  both sides of \eqref{h-1} and \eqref{h-2} satisfy the  same recurrences. We list these two recurrences as follows:
\begin{align*}
	&(2.13)^*: \quad  0=1-4(1+n)Sum[n]+2(5+3n)Sum[1+n]+(-5-2n)Sum[2+n],\\
	&(2.14)^*: \quad  0=1-2(1+n)Sum[n]+(4+3n)Sum[1+n]-(2+n)Sum[2+n].
\end{align*}	
It is routine to verify that both sides of \eqref{h-1} and \eqref{h-2} agree for all positive integers $n\geq 1$, respectively. We conclude that \eqref{h-1} and \eqref{h-2} holds for all $n\geq 1$. This concludes the proof.   \qed	
	
\section{Proof of the Theorem 1.1}	
According to the Dixon-Kummer summation theorem \cite{www-2-1},
\begin{align*}
{}_4F_3\left[
\begin{matrix}
	a,\ 1+\frac{a}{2},\ b,\ c \\
	\frac{a}{2},\ 1+a-b,\ 1+a-c
\end{matrix}; -1
\right]
=
\Gamma\left[
\begin{matrix}
	1+a-b,\ 1+a-c \\
	1+a,\ 1+a-b-c
\end{matrix}
\right],
\end{align*}	
Feng and Hou \cite{www-2}(or \cite{www-2-0}) found a WZ pair $(F(n,k), G(n,k))$ as follows,
\begin{align*}
F(n,k)=\frac{n(6n-4k-3)(2k)!(4n-2k-2)!(3n-k-2)!(2n-2k-1)!}{2^{6n-2k}(2n-k-1)!^2(n-k-1)!^2(2n)!(n-1)!k!},
\end{align*}	
\begin{align*}
	G(n,k)=-\frac{\alpha(n,k)(4n-2k)!(3n-k-1)!(2n-2k)!(2k)!}{2^{6n-2k+6}(2n-k+1)!(2n-k)!
		(n-k)^2(2n+1)!n!k!},
\end{align*}	
where 
\begin{align*}
	\begin{aligned}
\alpha(n,k)&=4k^4-44k^3n+180k^2n^2-308kn^3+184n^4-12k^3+98k^2n\\
&-260kn^2+214n^3+11k^2-62kn+79n^2-3k+9n,
\end{aligned}
\end{align*}
which satisfy the WZ equation
\begin{align}
	F(n+1,k)-F(n,k)=G(n,k+1)-G(n,k).  \label{03-a-1}
\end{align}

Our proof will start from equation \eqref{03-a-1}. Summing up \eqref{03-a-1} for $k$ from $0$ to $p-1$, we derive that
\begin{align}
	\sum_{k=0}^{p-1}F(n+1,k)-\sum_{k=0}^{p-1}F(n,k)=G(n,p)-G(n,0).  \label{03-a-2}
\end{align}
Since $G(n,k)=0$ for $n<k$ and $F(0,k)=0$, we continue to sum both sides of \eqref{03-a-2} over $n$ from $0$ to $\frac{p-3}{2}$, then	
\begin{align*}
	\sum_{k=0}^{p-1}F(\frac{p-1}{2},k)=-\sum_{n=0}^{p-1}G(n,0).
\end{align*}
	
Note that	
\begin{align*}
	&G(n,0)=-\frac{92n^2+61n+9}{192(2n+1)64^n}{2n \choose n}{3n \choose n}{4n \choose 2n},
\end{align*}
and	
\begin{align*}
	F(\frac{p-1}{2},k)=\frac{p-1}{8^{p-1}}{\frac{3p-5}{2}\choose p-1}\frac{4^k(3p-4k-6){2k\choose k}{2p-2k-4\choose p-k-2}{p-3-2k \choose \frac{p-3-2k}{2}}}{(3p-2k-5){\frac{3p-5}{2}\choose k}}.
\end{align*}	
Consequently, with some necessary arrangements, we obtain	
\begin{align}
	\begin{aligned}
		\sum_{k=0}^{\frac{p-3}{2}}&\frac{92k^2+61k+9}{(2k+1)64^k}{2k \choose k}{3k \choose k}{4k \choose 2k}\\
		&=\frac{192(p-1)}{8^{p-1}}{\frac{3p-5}{2} \choose p-1}\sum_{k=0}^{\frac{p-3}{2}}\frac{4^k(3p-4k-6){2k\choose k}{2p-2k-4 \choose p-k-2}{p-3-2k \choose \frac{p-3-k}{2}}}{(3p-2k-5){\frac{3p-5}{2}\choose k}}\\
		&=\frac{192(p-1)}{8^{p-1}}{\frac{3p-5}{2}\choose p-1}\mathcal{A}_1+\mathcal{A}_2, \label{03-1}
	\end{aligned}	
\end{align}		
where
\begin{align*}
	\mathcal{A}_1=\sum_{k=1}^{\frac{p-7}{2}}\frac{4^k(3p-4k-6){2k\choose k}{2p-2k-4\choose p-k-2}{p-3-2k\choose \frac{p-3-2k}{2}}}{(3p-2k-5){\frac{3p-5}{2}\choose k}},
\end{align*}	
\begin{align*}
		\mathcal{A}_2=\frac{72(p-1)^3}{8^{p-1}(2p-3)(3p-5)}{\frac{3p-5}{2}\choose p-1}{2p-2 \choose p-1}{p-1\choose \frac{p-1}{2}}
		+\frac{6p(p-1)}{4^{p-1}(p-2)}\left(1+\frac{(p-1)(p+4)}{(p+1)(p-4)}\right){p-1\choose \frac{p-1}{2}}^2.	
\end{align*}

Since
\begin{align}
	\begin{aligned}
		{\frac{3p-5}{2}\choose p-1}&=\frac{p(p+\frac{p-5}{2})(p+\frac{p-7}{2})\dots(p+1)(p-1)\dots(p-\frac{p+1}{2})}{(p-1)!}\\
		&\equiv \frac{(-1)^{\frac{p+1}{2}}2p(p+1)\left(1+pH_{\frac{p-5}{2}}+p^2H(1,1;\frac{p-5}{2})\right)\left(1-pH_{\frac{p+1}{2}}+p^2H(1,1;\frac{p+1}{2})\right)}{(p-1)(p-3){p-1 \choose \frac{p-1}{2}}}\\
		&\equiv \frac{(-1)^{\frac{p+1}{2}}2p\left(1+3p+9p^2-p^2H^{(2)}_{\frac{p-1}{2}}\right)}{3{p-1\choose \frac{p-1}{2}}}\\
		&\equiv -\frac{2p(1+3p+9p^2)}{4^{p-1}3} \pmod{p^4}. \label{03-2}
	\end{aligned}	
\end{align}
Therefore, we next dissect $\mathcal{A}_1$  modulo $p^3$ and $\mathcal{A}_2$ modulo $p^4$.

It is the well-know the congruence of Morley \cite{www-4}: for prime $p>2$,
\begin{align}
	\begin{aligned}
		{p-1\choose \frac{p-1}{2}}=(-1)^{\frac{p-1}{2}}4^{p-1} \pmod{p^3}. \label{03-3}
	\end{aligned}	
\end{align}
We also have
\begin{align}
	\begin{aligned}
		{2p-2\choose p-1}=\frac{p(p+1)\dots(p+p-2)}{(p-1)!}\equiv -p-2p^2 \pmod{p^3}. \label{a-03-3}
	\end{aligned}	
\end{align}
Combining  \eqref{03-3-1}, \eqref{03-2}, \eqref{03-3} and \eqref{a-03-3}, we deduce that
\begin{align}
	\begin{aligned}
\mathcal{A}_2\equiv -\frac{16p^2(15+49p)}{8^{p-1}\cdot75}\left(\frac{-1}{p}\right)+\frac{3}{8}4^{p-1}\left(16p-20p^2+11p^3\right) \pmod{p^4}. \label{03-4}
	\end{aligned}	
\end{align}

In view of \cite[Corollary 5.2]{www-1-4}:
\begin{align}
\sum_{k=1}^{(p-1)/2}\frac{(-1)^k}{k^2}\equiv 2(-1)^{\frac{p-1}{2}}E_{p-3} \pmod{p}.  \label{03-5}	
\end{align}
Furthermore, by necessary simplifying we obtain
\begin{align}
	\begin{aligned}
		\sum_{k=1}^{\frac{p-7}{2}}& \frac{(-1)^k(-76k^3-308k^2-387k-158)}{2(2k+1)^2(2k+3)^2(2k+5)^2(k+1)}\\
		&\equiv \frac{1}{12}\sum_{k=1}^{\frac{p-3}{2}}\frac{(-1)^k}{1+2k}+\frac{1}{6}H_{\frac{p-1}{2}}(-1) -\frac{7}{16}(-1)^{\frac{p-1}{2}}\sum_{k=1}^{\frac{p-1}{2}}\frac{(-1)^k}{k^2}-\frac{43}{576}(-1)^{\frac{p-1}{2}}+\frac{3157}{1800}\\		&\equiv \frac{3007}{1800}-\frac{q_p(2)}{6}+\frac{q_p(2)}{24}\left(\frac{-1}{p}\right)-\frac{43}{576}\left(\frac{-1}{p}\right)-\frac{7}{8}E_{p-3}  \pmod{p}, \label{03-6}
	\end{aligned}	
\end{align}	
where we used \eqref{0l-6}, \eqref{0l-7} and \eqref{03-5} in the last step.

For $1 \leq k \leq (p-3)/2$, we have
\begin{align*}
	\begin{aligned}	
		{2p-2k \choose p-k}
		&=\frac{p(2p-2k)(2p-2k-1)\dots(p+1)(p-1)\dots(p-k+1)}{(p-k)!}\\
		&\equiv {2k \choose k}\frac{(-1)^{k-1}p(1+pH_{p-2k})(1-pH_{k-1})}{{p-k \choose k}} \pmod{p^3}, 
	\end{aligned}	
\end{align*}
and
\begin{align*}
	\begin{aligned}	
		{p-k \choose k}&=\frac{(p-k)(p-k-1)\dots(p-2k+1)}{k!}\\
		&\equiv (-1)^k\frac{1}{2}{2k \choose k}(1-p(H_{2k-1}-H_{k-1})) \pmod{p^2}.
	\end{aligned}	
\end{align*}
It follows that
\begin{align}
	\begin{aligned}	
		{2k \choose k}{2p-2k \choose p-k} \equiv \frac{-2p(1-2pH_{k-1}+2pH_{2k-1})}{k} \pmod{p^3}. \label{03-7}
	\end{aligned}	
\end{align}
Meanwhile, for $1 \leq k \leq (p-7)/2$,
\begin{align}
	\begin{aligned}	
		{\frac{3p-5}{2} \choose k}&=\frac{(3p-1)(3p-1)\dots(3p-2k-3)}{\left(\frac{3p}{2}-(k+1)\right)\left(\frac{3p}{2}-k\right)\dots(\frac{3p}{2}-1)k!(3p-1)(3p-3)2^{2k+1}}\\
		&\equiv   \frac{(-1)^k(2k+1)(2k+3){2k\choose k}}{2^{2k}(3p-1)(3p-3)\left(1-\frac{3}{2}pH_{k+1}(1)+3pH_{2k+3}(1)\right)} \pmod{p^2}. \label{03-8}
	\end{aligned}	
\end{align}
Similarly, it is easy for us to check that
\begin{align*}	
	{p-1\choose 2k}&\equiv (-1)^k\left(1-pH_{2k}\right) \pmod{p^2},
\end{align*}
and
\begin{align*}	
	{\frac{p-1}{2}\choose k}&\equiv {2k\choose k}\left(1-p\sum_{i=1}^{k}\frac{1}{2i-1}\right)/(-4)^k \pmod{p^2}.
\end{align*}
Combining \eqref{03-3} and the above congruence gives
\begin{align}
	\begin{aligned}	
	{p-1-2k\choose \frac{p-1-2k}{2}}&=\frac{{p-1\choose \frac{p-1}{2}}{\frac{p-1}{2}\choose k}^2}{{2k\choose k}{p-1\choose 2k}}
	\equiv \frac{(-1)^{\frac{p-1}{2}}4^{p-1}{2k\choose k}}{16^k}\left(1-pH_{2k}(1)+pH_k(1)\right) \pmod{p^2}. \label{03-9}
\end{aligned}
\end{align}

Then, using \eqref{03-7}, \eqref{03-8} and \eqref{03-9} together with simplification, we obtain
\begin{align*}
	\begin{aligned}
		\mathcal{A}_1&\equiv  \frac{3}{8}(-1)^{\frac{p-1}{2}}4^{p-1}\sum_{k=1}^{\frac{p-7}{2}} \frac{(-1)^kp(4p-1)}{(2k+1)(2k+3)(2k+5)}\\
		&+ \frac{3}{8}(-1)^{\frac{p-1}{2}}4^{p-1}p^2\sum_{k=1}^{\frac{p-7}{2}} \frac{(-1)^k(-76k^3-308k^2-387k-158)}{2(2k+1)^2(2k+3)^2(2k+5)^2(k+1)}     \\
		&+\frac{3}{8}(-1)^{\frac{p-1}{2}}4^{p-1}p^2\left(\sum_{k=1}^{\frac{p-7}{2}} \frac{5(-1)^kH_k}{2(2k+1)(2k+3)(2k+5)}-\sum_{k=1}^{\frac{p-7}{2}} \frac{4(-1)^kH_{2k}}{(2k+1)(2k+3)(2k+5)}\right) \pmod{p^3}.
	\end{aligned}	
\end{align*}		
Substituting \eqref{0l-1-0}, \eqref{0l-1}, \eqref{0l-2}, \eqref{03-6} into the above expression and simplifying yields
\begin{align}
	\begin{aligned}
		\mathcal{A}_1
		&\equiv \frac{3}{8}4^{p-1}\left((\frac{41}{32}p^2-\frac{1}{4}p)q_p(2)-\frac{7}{8}p^2q_p(2)^2+(\frac{2}{5}p+\frac{451}{900}p^2-\frac{p^2}{4}E_{p-3})\left(\frac{-1}{p}\right)\right)\\
		&-\frac{3}{8}4^{p-1}(\frac{51}{128}p^2-\frac{5}{32}p) \pmod{p^3}.
	\end{aligned}	
\end{align}	
	
For  any positive integer $a$, we have 
\begin{align*}
(2^a)^{p-1}=\left(1+q_p(2)p\right)^a. 
\end{align*}
Finally, substituting \eqref{03-2}, \eqref{03-4} and \eqref{03-9} into \eqref{03-1} and applying Fermat's Little Theorem, we arrive at \eqref{03} after concrete caculations. This concludes the proof. \qed

\section{Proof of the Theorem 1.2}	
	
Summing both sides of \eqref{03-a-1} over $n$ from $0$ to $p-1$ yields
\begin{align*}
	F(p,k)-F(0,k)=\sum_{n=0}^{p-1}G(p,k+1)-\sum_{n=0}^{p-1}G(p,k).
\end{align*}	
Continuing to sum over $k$ from $0$ to $p-1$, utilizing  $G(n,k)=0$ for $n<k$ and $F(0,k)=0$, we obtain	
\begin{align*}
	\sum_{k=0}^{p-1}F(n,k)=-\sum_{n=0}^{p-1}G(p,k),
\end{align*}		
which implies
\begin{align}
	\sum_{k=0}^{p-1}\frac{92k^2+61k+9}{64^k(2k+1)}{2k \choose k}{3k\choose k}{4k\choose 2k}=\frac{192p{3p-1\choose 2p}}{64^p}\mathcal{B}, \label{04-0}
\end{align}	
where 	
\begin{align*}
	\mathcal{B}=\sum_{k=0}^{p-1}\frac{(6p-4k-3)(p-k)4^k{2k \choose k}{2p-2k\choose p-k}{4p-2k-2\choose 2p-k-1}}{2(3p-k-1)(2p-2k-1){3p-1\choose k}}.
\end{align*}		
	
For $(p+1)/2\leq k \leq p-2$, we have
\begin{align}
	\begin{aligned}
	{2k\choose k}{2p-2k\choose p-k}&={2k\choose k}\frac{(2p-2k)(2p-2k-1)\dots(2p-(p+k-1))}{(p-k)!}\\
	&\equiv \frac{2p}{k}\left(1-2p\sum_{i=2k}^{p+k-1}\frac{1}{i}\right) \pmod{p^3}. \label{04-1}
\end{aligned}	
\end{align}
And for $1\leq k \leq p-2$, it is easy to see that
\begin{align}
	\begin{aligned}
		{4p-2k-2 \choose 2p-k-1}&=\frac{2p(2p+2p-k-1)\dots(2p+1)(2p-1)\dots(2p-k)}{(2p-k-1)!}\\
		&\equiv \frac{2(-1)^kp}{(k+1){2p-k-1\choose k+1}} \pmod{p^2}. \label{04-2}
	\end{aligned}	
\end{align}
Meanwhile, for $1\leq k \leq p-2$ and $k\neq \frac{p-1}{2}$, the polynomial $(3p-1-k)(2p-2k-1){3p-1\choose k}$ does not contain $p$ as a factor. Hence, combining \eqref{03-7}, \eqref{04-1} and \eqref{04-2} gives
\begin{align}
	\begin{aligned}
	\mathcal{B}&\equiv \frac{(6p-3){4p-2\choose p-1}{2p-2\choose p-1}}{3p-1}+\frac{2(4p-1)4^{p-1}{p-1\choose (p-1)/2}^2{3p-1\choose (3p-1)/2}}{(5p-1){3p-1\choose (p-1)/2}}\\
	&+\frac{4^{p-1}(4p^2-1){2(p-1)\choose p-1}^2}{p^2{3p-1\choose p-1}} \pmod{p^2}.  \label{04-4}
	\end{aligned}
\end{align}

Note that
\begin{align*}
	\begin{aligned}
	{3p-2\choose \frac{3p-1}{2}}&=2\frac{(2p+p-1)\dots(2p+1)(p+p-1)\dots(p+\frac{p+1}{2})}{(p+\frac{p-1}{2})\dots(p+1)(p-1)!}\\
	&\equiv 2\left(1+2p\sum_{k=\frac{p+1}{2}}^{p-1}\frac{1}{j}\right){p-1\choose \frac{p-1}{2}} \pmod{p^2},
	\end{aligned}
\end{align*}
similarly, we deduce that
\begin{align*}
	\begin{aligned}
	{3p-1\choose (3p-1)/2}/{3p-1\choose (p-1)/2}\equiv 2\left(1+pH_{p-1}\right) \pmod{p^2}.
	\end{aligned}
\end{align*}	
Then, utilizing \eqref{03-3} and \eqref{a-03-3}, we deduce that
\begin{align}
	\begin{aligned}
	2(4p-1)4^{p-1}{p-1\choose (p-1)/2}^2{3p-1\choose (3p-1)/2}/\left((5p-1){3p-1\choose (p-1)/2}\right)    \equiv 4(p+1)32^{p-1} \pmod{p^2},  \label{04-5}
	\end{aligned}
\end{align}
and
\begin{align}
	\begin{aligned}
	4^{p-1}(4p^2-1){2(p-1)\choose p-1}^2/\left(p^2{3p-1\choose p-1}\right)\equiv -4^{p-1}\left(1+4p\right)    \pmod{p^2}.  \label{04-6}
	\end{aligned}
\end{align}
It is easy to verify that
\begin{align}
	\begin{aligned}
	(6p-3){4p-2\choose p-1}{2p-2\choose p-1}/(3p-1)\equiv 0 \pmod{p^2}.  \label{04-7}
	\end{aligned}
\end{align}	

Substituting \eqref{04-5}--\eqref{04-7}	into \eqref{04-4} and simplifying, we obtain
\begin{align*}
	\frac{192p{3p-1\choose 2p}}{64^p}\mathcal{B}\equiv \frac{12(p+p^2)}{2^{p-1}}-\frac{3(p+4p^2)}{16^{p-1}}\equiv 9p^2 \pmod{p^3}.
\end{align*}	
This, together with \eqref{04-0} gives \eqref{04}. 

Thus, we complete the proof Theorem $1.2$. \qed

\section{Proof of the Corollary 1.3}	
On the one hand, note that for $1 \leq k \leq \frac{p-1}{2}$, we have
\begin{align*}
	\begin{aligned}
	&{2p+2k \choose p+k}/{3p-1 \choose p-1-k}\\
	&=\frac{2p(2p+1)\dots(2p+1)(2p-1)\dots(2p-(p-k-1))(2p+p-1)\dots(2p+k+1)}{(p+k)!(p-1-k)!}
	\\
	&\equiv \frac{2p(-1)^k}{(2k+1){p-1\choose k}{p+k\choose 2k+1}}\left(1-2p\sum_{i=k+1}^{p-1}\frac{1}{i}+2p\sum_{i=1}^{2k}\frac{1}{i}-2p\sum_{k=1}^{p-1-k}\frac{1}{i}\right)\\
	&\equiv \frac{2p(-1)^k}{(2k+1){p-1\choose k}{p+k\choose 2k+1}}\left(1+2p\sum_{i=1}^{2k}\frac{1}{i}\right) \pmod{p^2}.
	\end{aligned}	
\end{align*}
This, with together \eqref{03-7} yields
\begin{align*}
	\begin{aligned}
		\sum_{k=1}^{\frac{p-3}{2}}&\frac{(2p+4k+1)(p-k)4^{p-1-k}{2k \choose k}{2p-2k\choose p-k}{2p+2k\choose p+k}}{2(2p+k)(2p-2k-1){3p-1\choose p(n)-1-k}}
		\equiv p\sum_{k=1}^{\frac{p-3}{2}}\frac{2p(-1)^k(4k+1)}{(2k+1)k^2{p-1\choose k}{p-1+k\choose 2k}} \pmod{p^2},
	\end{aligned}	
\end{align*}
where we also used Fermat's Little Theorem.	Moreover, it is easy to see that
\begin{align*}
	\begin{aligned}
		{p-1+k\choose 2k}\equiv \frac{(-1)^kp}{k{2k\choose k}} \pmod{p^2} \quad and \quad {p-1\choose k}\equiv (-1)^k\pmod{p}.
	\end{aligned}	
\end{align*}		
It follows that
\begin{align}
	\begin{aligned}
		\sum_{k=1}^{\frac{p-3}{2}}&\frac{(2p+4k+1)(p-k)4^{p-1-k}{2k \choose k}{2p-2k\choose p-k}{2p+2k\choose p+k}}{2(2p+k)(2p-2k-1){3p-1\choose p(n)-1-k}}\\
		&\equiv 2p\sum_{k=1}^{\frac{p-3}{2}}\frac{{2k\choose k}(1+4k)}{(-4)^kk(1+2k)}  \pmod{p^2}. \label{05-1}
	\end{aligned}	
\end{align}	
	
On the other hand, by \eqref{04-1} and \eqref{04-2} we have	
\begin{align}
	\sum_{k=\frac{p+1}{2}}^{p-2}\frac{(6p-4k-3)(p-k)4^k{2k \choose k}{2p-2k\choose p-k}{4p-2k-2\choose 2p-k-1}}{2(3p-k-1)(2p-2k-1){3p-1\choose k}}\equiv 0 \pmod{p^2}. \label{05-2}
\end{align}	
We replace $k$ with $p-1-k$ on the left-side of \eqref{05-2}, then	
\begin{align}
	\begin{aligned}
		\sum_{k=\frac{p+1}{2}}^{p-2}&\frac{(6p-4k-3)(p-k)4^k{2k \choose k}{2p-2k\choose p-k}{4p-2k-2\choose 2p-k-1}}{2(3p-k-1)(2p-2k-1){3p-1\choose k}}\\
		&=\sum_{k=1}^{\frac{p-3}{2}}\frac{(2p+4k+1)(p-k)4^{p-1-k}{2k \choose k}{2p-2k\choose p-k}{2p+2k\choose p+k}}{2(2p+k)(2p-2k-1){3p-1\choose p(n)-1-k}}. \label{05-3}
	\end{aligned}	
\end{align}	
According to \eqref{05-1}--\eqref{05-3} we get
\begin{align}
\sum_{k=1}^{\frac{p-3}{2}}\frac{{2k\choose k}(1+4k)}{(-4)^kk(1+2k)} \equiv 0 \pmod{p}. \label{05-4}	
\end{align}		
	
A result due to Sun \cite[(4.4)]{www-1-2} says: for $p=2n+1$ be an odd prime,	
\begin{align}
	\frac{(-4)^k}{{2k\choose k}}{n\choose k}\equiv 1-p\sum_{j=1}^{k}\frac{1}{2j-1} \pmod{p^2}. \label{05-5}	
\end{align}	
Therefore, combining \eqref{0l-4},  \eqref{h-1}, \eqref{h-2}, \eqref{03-3}, \eqref{05-4} and \eqref{05-5} we get
\begin{align*}
	\begin{aligned}
	\sum_{k=1}^{\frac{p-3}{2}}\left(\frac{1}{k}+\frac{2}{1+2k}\right){\frac{p-1}{2}\choose k}&=\sum_{k=1}^{\frac{p-1}{2}}\frac{2^k}{k}+\frac{2^{p}}{p{p-1\choose \frac{p-1}{2}}} \left(\sum_{k=1}^{\frac{p-1}{2}}\frac{{2k\choose k}}{2^k}+1\right)-H_{\frac{p-1}{2}}
	-\frac{2(p^2+p-1)}{p(p-1)}\\
	&\equiv \frac{2(-1)^{\frac{p-1}{2}}}{2^{p-1}p}\sum_{k=1}^{\frac{p-1}{2}}\frac{{2k\choose k}}{2^k}+\sum_{k=1}^{\frac{p-1}{2}}\frac{2^k}{k}+\frac{2(-1)^{\frac{p-1}{2}}-2^p}{2^{p-1}p}+2q_p(2)\\
	&\equiv 0  \pmod{p}.
\end{aligned}
\end{align*}
From the above congruence, after explicit calculations, we obtain
\begin{align*}
	\sum_{k=1}^{\frac{p-1}{2}}\frac{{2k \choose k}}{2^k}\equiv  -\frac{p}{2}\left(\frac{-1}{p}\right)\sum_{k=1}^{\frac{p-1}{2}}\frac{2^k}{k}+\left(\frac{-1}{p}\right)-1 \pmod{p^2}, 	
\end{align*}
which also  means
\begin{align*}
	\sum_{k=1}^{\frac{p-1}{2}}\frac{{2k \choose k}}{2^k}\equiv  \left(\frac{-1}{p}\right)-1 \pmod{p}.	
\end{align*}

The proof of Corollary $1.3$ is now complete. \qed

\end{document}